\documentclass[a4paper,twoside,11pt,reqno]{amsart}
\newcommand{\Ueberschrift}{Birational $p$-adic Galois sections in higher dimensions}
\newcommand{\Kurztitel}{Birational $p$-adic Galois sections}
\usepackage{amsmath} \usepackage{amssymb} \usepackage{amsopn}
\usepackage{amsthm} \usepackage{amsfonts} \usepackage{mathrsfs}

\usepackage[headings]{fullpage}
\usepackage[bottom]{footmisc}

\usepackage[dvips]{graphicx} \usepackage[usenames]{color}
\usepackage[all]{xy}
\usepackage[OT2, T1]{fontenc}
\usepackage[utf8]{inputenc}
\usepackage{dsfont}

\usepackage[pdfpagelabels, pdftex]{hyperref}

\hypersetup{
  pdftitle={},
  pdfauthor={},
  pdfsubject={},
  pdfkeywords={},
  colorlinks=true,
  breaklinks=true,
  bookmarksopen=true,
  bookmarksnumbered=true,
  pdfpagemode=UseOutlines,
  plainpages=false}

\newcommand{\bP}{{\mathbb P}}
\newcommand{\bQ}{{\mathbb Q}}

\newcommand{\bZ}{{\mathbb Z}}

\newcommand{\cS}{{\mathscr S}}

\newcommand{\dO}{{\mathcal O}}

\newcommand{\fp}{{\mathfrak p}}

\DeclareSymbolFont{cyrletters}{OT2}{wncyr}{m}{n}
\DeclareMathSymbol{\Sha}{\mathalpha}{cyrletters}{"58}

\newcommand{\surj}{\twoheadrightarrow} 
\newcommand{\inj}{\hookrightarrow}

\DeclareMathOperator{\coker}{coker}
\DeclareMathOperator{\im}{im}

\DeclareMathOperator{\Spec}{Spec}

\newcommand{\OO}{\dO}

\DeclareMathOperator{\Alb}{Alb}

\DeclareMathOperator{\res}{res}

\DeclareMathOperator{\Gal}{Gal}

\newcommand{\ab}{{\rm ab}}

\newcommand{\abs}{{\rm abs}}

\newcommand{\ov}[1]{\mbox{${\overline{#1}}$}} 

\newtheorem{thm}{Theorem}

\newtheorem{prop}[thm]{Proposition}
\newtheorem{lem}[thm]{Lemma}

\newtheorem*{thmMAIN}{Theorem}

\theoremstyle{definition}

\theoremstyle{remark}
\newtheorem{rmk}[thm]{Remark}

\newtheorem{ex}[thm]{Example}

\newenvironment{pro}[1][Proof]{{\it{#1:}} }{\hfill $\square$}
\newenvironment{pro*}[1][Proof]{{\it{#1:}} }{}
\newenvironment{pro**}[1][]{{\it{#1}} }{\hfill $\square$}

\numberwithin{equation}{section}

\begin{document}

\hrule width\hsize
\hrule width\hsize
\hrule width\hsize
\hrule width\hsize

\vspace{1cm}

\title[\Kurztitel]{\Ueberschrift} 
\author{Jakob Stix}
\address{Jakob Stix, MATCH - Mathematisches  Institut, Universit\"at Heidelberg, Im Neuenheimer Feld 288, 69120 Heidelberg, Germany}
\email{stix@mathi.uni-heidelberg.de}
 \urladdr{http://www.mathi.uni-heidelberg.de/~stix/}

\date{\today} 

\maketitle

\begin{quotation} 
  \noindent \small {\bf Abstract} ---  This note explores the consequences of Koenigsmann's model theoretic argument from the proof of the birational $p$-adic section conjecture for curves in the context of higher dimensional varieties over $p$-adic local fields. 
  \end{quotation}

%%%%%%%%%%%%%%%%%%%%%%%%%%%%%%%%%%%%%%%%%%%
%%%%%% Main body                 %%%%%%%%%%%%%%%%%%%%%%%%%%%
%%%%%%%%%%%%%%%%%%%%%%%%%%%%%%%%%%%%%%%%%%%

%%%%%%%%%%%%%%%%%%%%%%%%%%%%%%%%%%%%%%%%%%
\section{Introduction}
%%%%%%%%%%%%%%%%%%%%%%%%%%%%%%%%%%%%%%%%%%

A \textbf{birational Galois section}  for a geometrically irreducible and reduced variety $X/k$ is a continuous section of the restriction homomorphism 
\[
\res_{K/k} \, : \ \Gal_K \to \Gal_k,
\] 
where $\Gal_K$ and $\Gal_k$ are the absolute Galois groups of the function field $K = k(X)$ and of $k$, with respect to an algebraic closure $\ov{K}$ of $K$ and the algebraic closure $\bar{k}$ of $k$ contained in $\ov{K}$. 

A proof of Grothendieck's anabelian section conjecture would imply, for a smooth curve $X/k$ over a finitely generated extension $k/\bQ$, that birational Galois sections come in packets indexed by the $k$-rational points $X(k)$. In particular, then $X(k) \not= \emptyset$ 
%has a $k$-rational point 
if and only if $X/k$ admits birational Galois sections. These consequences of the section conjecture is known in some cases for number fields $k$ by Stoll \cite{stoll:finitedescent} Remark 8.9, see also \cite{hararistix:findesc} Theorem 17.

More strikingly, in \cite{koenigsmann:birationalsc} Proposition 2.4 (b), Koenigsmann proves, with a touch of model theory of $p$-adically closed fields, that birational Galois sections for smooth curves $X/k$ over finite extensions $k/\bQ_p$ come in packets indexed by $X(k)$. The packet of $a \in X(k)$ consists of those sections $s$ with $s(\Gal_k) \subset D_{\bar v}$, where $D_{\bar v}$ is any  decomposition group of a prolongation $\bar{v}$ to $\ov{K}$  of the $k$-valuation $v$ corresponding to $a$ with valuation ring $\OO_{X,a} \subset k(X)$.

\smallskip

Our main result is the following theorem, see Theorem \ref{thm:main}.

\begin{thmMAIN}
Let $X/k$ be a geometrically irreducible, normal variety over a finite extension $k/\bQ_p$ with function field $K$. Then every  birational Galois section 
has image in the decomposition subgroup $D_{\bar v} \subset \Gal_K$ for a unique $k$-valuation
 $\bar{v}$ of $\ov{K}$ with residue field of $v = \bar{v}|_K$ equal to  $k$. 
 
In particular, conjugacy classes of sections of $\Gal_K \to \Gal_{k}$ come in disjoint packets associated to $k$-valuations $v$ of $K$ with residue field $k$.
\end{thmMAIN}

%%%%%%%%%%%%%%%%%%%%%%%%%%%%%%%%%%%%%%%%%%%
%%%%%%%%%%%%%%%%%%%%%%%%%%%%%%%%%%%%%%%%%%%

\medskip

\noindent
{\bf Acknowledgements.}
I thank Zo\'e Chatzidakis and Fran\c{c}ois Loeser for inviting me to speak in the seminar \textit{G\'eom\'etrie et Th\'eorie des Mod\`eles} and I am grateful to Jean-Louis Colliot-Th\'el\`ene and Jochen Koenigsmann  for a subsequent discussion on Koenigsmann's Lemma.

%%%%%%%%%%%%%%%%%%%%%%%%%%%%%%%%%%%%%%%%%%%
%%%%%%%%%%%%%%%%%%%%%%%%%%%%%%%%%%%%%%%%%%%

%%%%%%%%%%%%%%%%%%%%%%%%%%%%%%%%%%%%%%%%%%%
\section{Review of Koenigsmann's use of model theory}
%%%%%%%%%%%%%%%%%%%%%%%%%%%%%%%%%%%%%%%%%%%

We keep the notation from the introduction and refer to \cite{prestelroquette} for the notion of a $p$-adically closed field. 
Koenigsmann has the following lemma in the case of smooth curves, \cite{koenigsmann:birationalsc} Proposition 2.4 (a). The general case, see  \cite{koenigsmann:birationalsc} Remark 2.5, admits a mathematically identical proof. We decide to nevertheless give a proof in order to hopefully make the argument more transparent for non-model theorists.

\begin{lem}[Koenigsmann's Lemma]
Let $k$ be a $p$-adically closed field, and let $X/k$ be a geo\-metri\-cally irreducible and reduced variety. If $X/k$ admits a birational Galois section, then $X$ has a $k$-rational point.
\end{lem}
\begin{pro}
Let $s: \Gal_k \to \Gal_K$ be a section. We set $L= \ov{K}^{s(\Gal_k)}$ for the fixed field of  the image, so that by construction, the restriction map $\Gal_L \to \Gal_k$ is an isomorphism. Thus  $k \subset L$ is relatively algebraically closed and the subfields 
$k^{\abs} = L^{\abs}$ of absolutely algebraic elements, i.e., algebraic over the prime field, agree. Since the restriction map $\Gal_k \to \Gal_{k^\abs}$ is an isomorphism \cite{pop:galoiskennzeichnung} (E4), we deduce that $\Gal_L \to \Gal_{L^\abs}$ is an isomorphism as well. Now \cite{pop:galoiskennzeichnung} Theorem (E12) implies that $L$ is $p$-adically closed with the same invariants $e$ and $f$ as $k$. 

\smallskip
 
With Prestel and Roquette \cite{prestelroquette} Theorem 5.1 we conclude that $k \subset L$ is an elementary extension with respect to the model theory of valued fields. This implies that the statement in the language of fields with constants in $k$ saying 
\[
\text{\textit{ 'the set of rational points of $X$ is non-empty'}}
\]
is true over $k$ if and only if it is true over $L$ where the constants from $k$ are interpreted via the embedding $k \subset L$. Since $X$ admits the tautological $L$-rational point 
\[
\Spec(L) \to \Spec(K) \to X,
\]
we are done.
\end{pro}

%%%%%%%%%%%%%%%%%%%%%%%%%%%%%%%%%%%%%%%%%%%
\section{From existence of points to control of unramified quotients}
%%%%%%%%%%%%%%%%%%%%%%%%%%%%%%%%%%%%%%%%%%%

The fundamental group of $X/k$ fits naturally in an extension $\pi_1(X/k)$
\[
1 \to \pi_1(X_{\bar k}) \to \pi_1(X) \to \Gal_k \to 1,
\]
where  the geometric generic point $\Spec(\ov{K}) \to X_{\bar k} \to X$  is the implicit base point. The space of sections of $\pi_1(X/k)$ up to conjugation by elements from $\pi_1(X_{\bar k})$ will be denoted by 
\[
\cS_{\pi_1(X/k)}
\]
and its birational analogue, the set of $\Gal_{K\bar{k}}$-conjugacy classes of sections of $\Gal_K \to \Gal_k$, will be denoted by 
\[
\cS_{\pi_1(K/k)}.
\]
By functoriality, to $a \in X(k)$ we associate a class of sections $s_a : \Gal_k \to \pi_1(X)$. This gives rise to the non-abelian Kummer map $a \mapsto \kappa(a)  = s_a$
\[
\kappa  \, : \  X(k) \to \cS_{\pi_1(X/k)}
\]
Let $j: \Spec K \to X$ be the inclusion of the generic point which on $\pi_1$ induces the  map 
\[
j_\ast : \Gal_K \to \pi_1(X),
\]
a surjection if $X$ is normal, 
and furthermore a map
\[
j_\ast : \cS_{\pi_1(K/k)} \to \cS_{\pi_1(X/k)},
\]
the image of which by definition is the set of \textbf{birationally liftable} sections.

\begin{prop} \label{prop:imkappa=imj}
Let $k$ be a finite extension of $\bQ_p$, and let $X/k$ be a proper, smooth and geometrically irreducible variety. Then the image of the non-abelian Kummer map 
\[
\kappa \, : \  X(k) \to \cS_{\pi_1(X/k)}
\]
agrees with the image of the natural map 
\[
j_\ast : \cS_{\pi_1(K/k)} \to \cS_{\pi_1(X/k)}.
\]
\end{prop}
\begin{pro}
Recall the notion of a neighbourhood of a section $s: \Gal_k \to \pi_1(X)$. This is a finite \'etale cover $X' \to X$ together with a lift $s'$ of the section, i.e., an open subgroup $H \subseteq \pi_1(X)$ containing the image of the section $s=s'$. The limit over all neighbourhoods yields a pro-\'etale cover 
\[
X_s = \varprojlim X' \to X
\]
corresponding to $\pi_1(X_s) = s(\Gal_k) \subseteq \pi_1(X)$. It follows that a section $s$ of $\pi_1(X/k)$ comes from a $k$-rational point if and only if $X_s(k)$ is nonempty. More precisely, we have $s = s_a$ for every 
\[
a \in \im\big(X_s(k) \to X(k)\big).
\]
We have 
\[
X_s(k) = \varprojlim X'(k)
\]
where $X'$ ranges over all neighbourhoods of the section $s$. Since $X/k$ is assumed proper, all sets $X'(k)$ are naturally compact in the $p$-adic topology. Thus the projective limit is nonempty if and only if it is a projective limit of nonempty sets.  We conclude that  $s = s_a$ for some $a \in X(k)$ if and only if all neighbourhoods $(X',s')$ of $s$ have $k$-rational points. The latter follows from being in the image of $j_\ast$ by Koenigsmann's Lemma above because with $s$ also $s'$ is birationally lifting as a section of $\pi_1(X'/k)$.

\smallskip

To conclude the converse, it suffices to construct a $k$-valuation\footnote{ A $k$-valuation is a valuation that restricts to the trivial valuation on $k$, i.e., with $k$ in the valuation ring.}
 $v$ on $K=k(X)$ such that 
\begin{enumerate}
\item[(i)]
a given point $a \in X(k)$ is the center of the valuation on $X$,
\item[(ii)]
the residue field of $v$ agrees with $k$,
\item[(iii)] 
and the natural surjection $D_v \surj \Gal_k$ splits where  $D_v \subset \Gal_K$ is the  decomposition subgroup of the valuation.
\end{enumerate}
Then the decomposition subgroup projects to $j_\ast(D_v) = s_a(\Gal_k)$, and thus any splitting of $D_v \to \Gal_k$ will allow to lift $s_a$ to a section of $\Gal_K \to \Gal_k$.

We construct such a valuation $v$. Let $t_1, \ldots, t_d$ be a system of parameters in the regular local ring $\OO_{X,a}$ and put $\fp_i = (t_1,\ldots, t_i)$ for $0 \leq i \leq d$ corresponding to generic points $a_i$ of irreducible cycles $Z_i \subseteq X$ which are regular at $a=a_d$. We consider the valuation $v$ associated to the corresponding Parshin chain 
\[
X = Z_0 \supset Z_1 \supset \ldots \supset Z_d = a
\]
which is the composition in the sense of valuations of the discrete rank $1$ valuation rings
$\OO_{Z_i,a_{i+1}}$.
Its decomposition group $D_v$ sits in an extension
\[
1 \to \hat{\bZ}(1)^d \to D_v \to \Gal_k \to 1 
\]
which is split by choosing compatible $n$th roots for all parameters $t_i$ and all $n$. 
\end{pro}

\begin{rmk}
(1) 
For birationally liftable sections to be contained in $\kappa(X(k))$ requires $X/k$ to be proper, or at least a smooth compactification with no $k$-rational point in the boundary. Any affine smooth hyperbolic curve with a $k$-rational point  at infinity provides a counter example.

(2)
That sections $s_a$ are birationally liftable requires a certain amount of local regularity at $k$-rational points. It is not enough to assume $X$ normal as the following example shows.
\end{rmk}

\begin{ex}
Let $k/\bQ_p$ be finite and let $C/k$ be a smooth projective geometrically irreducible curve which is 
\begin{enumerate}
\item[(i)]
not hyperelliptic,
\item[(ii)]
and has no $k$-rational point.
\end{enumerate}
We consider the difference map  $(a,b) \mapsto a-b$ to the Albanese variety $\Alb_C$ of $C$.
\[
d: C \times C \to A = \Alb_C
\]
We set $f : C \times C \to X$ for the Stein-factorization of $d$. This $X$ is a normal variety and, due to assumption (i), the map $d$  contracts\footnote{ This $\Delta$ has negative self-intersection and thus can be contracted to an algebraic space anyway; but we would like to have a normal algebraic variety.}
 exactly the diagonally embedded 
\[
\Delta = C \inj C \times C.
\]
So $X$ can be identified with the contraction of $\Delta$ to a $k$-rational point $\star$
\[
X = C \times C/\Delta \sim \star,
\]
thus $X$ is birational to $C \times C$.
If $X/k$ were eligible for the conclusion of Proposition \ref{prop:imkappa=imj}, then the section $s_\star$ associated to $\star \in X(k)$ admits birational lifting and in particular lifts to a birationally lifting section $s$ of $\pi_1(C \times C /k)$. By Proposition \ref{prop:imkappa=imj} applied to the smooth projective $C \times C$ we find a $k$-rational point $(a,b) \in C\times C(k)$ which contradicts assumption (ii). 

\medskip

We now forget about assumption (ii) but keep (i). As an aside, we can determine the fundamental group extension $\pi_1(X/k)$ as follows. 
The induced maps on fundamental groups are surjective 
\[
\pi_1(C \times C) \surj \pi_1(X) \surj \pi_1(A).
\]
Since the diagonally embedded copy of $\pi_1(C_{\bar k})$ goes to zero, the computation in $\ker(\pi_1(f))$
\[
(1,h)(g,g)(1,h)^{-1} (g,g)^{-1} = (1,hgh^{-1}g^{-1})
\]
shows that we more precisely have 
\[
\coker\big(\pi_1(\Delta) :  \pi_1^\ab(C_{\bar k}) \to \pi^\ab_1(C_{\bar k}) \times \pi_1^\ab(C_{\bar k}) \big) \surj  \pi_1(X_{\bar k}) \surj \pi_1(A_{\bar k})
\]
with the composite map being an isomorphism (by geometric class field theory). In particular, the map $X \to A$ induces a canonical isomorphism
\[
\pi_1(X/k) = \pi_1(A/k)
\]
as extensions. The image of  the $p$-adic analytic map 
\[
d \, : \ C(k) \times C(k) \to  X(k) \to A(k) 
\]
lies in a $2$-dimensional $p$-adic analytic subspace of the $g \geq 3$-dimensional $p$-adic manifold $A(k)$. Therefore, the map $d$ is never surjective on $k$-rational points. Sections $s_a$ for $\pi_1(A/k)$ for $k$-rational points in $a \in A(k) \setminus X(k)$ are thus sections of $\pi_1(X/k)$ which are not birationally liftable and for which there is no $k$-rational point on $X$ that is responsible for it.
\end{ex}

%%%%%%%%%%%%%%%%%%%%%%%%%%%%%%%%%%%%%%%%%%%
\section{On birational sections and valuations}
%%%%%%%%%%%%%%%%%%%%%%%%%%%%%%%%%%%%%%%%%%%

We keep the notation from above, but now assume $X$ is normal. A neighbourhood of a section $s : \Gal_k \to \Gal_K$ is given by a normal finite branched cover $X' \to X$, i.e., a finite extension $K'/K$ inside $\ov{K}$, such that the image of $s$ lies in $\Gal_{K'} \subseteq \Gal_K$. The limit over all these neighbourhoods yields a pro-branched cover again denoted
\[
X_s = \varprojlim X' \to X
\]
corresponding to $\pi_1(k(X_s)) = s(\Gal_k) \subset \Gal_K$. 

\begin{prop} \label{prop:card1}
Let $k$ be a finite extension of $\bQ_p$, and let $X/k$ be a normal projective, geometrically irreducible variety. Then the set $X_s(k)$ consists of exactly one element.
\end{prop}
\begin{pro}
If we can show that the image $X_s(k) \to X(k)$ consists of only one element, then we can apply this to all neighbourhoods $X'$ of $s$ and achieve the proof of the proposition. 

That the image is nonempty follows from one half of the proof of  Proposition  \ref{prop:imkappa=imj} above, namely the half which only requires $X/k$ to be proper.

A point $a \in \im\big(X_s(k) \to X(k)\big)$ necessarily has $j_\ast(s) = s_a$ as sections of $\pi_1(X/k)$. 
So let us assume that we have $a \not= b$ in the image. By the projective version of Noether's Normalization Lemma, we can choose a suitable finite map 
\[
f : X \to \bP_k^n
\]
with $n = \dim(X)$, namely 
 by first choosing an immersion $X \inj \bP^N_k$ and then doing a sufficiently generic linear projection 
(note that $k$ is an infinite field to facilitate matters). 
We furthermore assume that $f(a) \not= f(b)$, which is an open condition. Next, we pick our preferred abelian variety $A/k$ of dimension $n$ and apply Noether's Normalization Lemma again to obtain a finite map 
\[
A \to \bP^n_k.
\]
Let $X'$ be defined as the normalization of $X$ in the compositum of function fields of $X$ and of $A$ over that of $\bP_k^n$, i.e., we have a commutative (not necessarily cartesian) square
\[
\xymatrix@M+1ex@R-2ex{
X' \ar[r]^{f'} \ar[d] & A \ar[d] \\
X \ar[r]^f & \bP_k^n.
}
\]
The branched cover $X' \to X$ however, as can be seen by looking at $\Gal_{K'} \subseteq \Gal_K$ where $K'$ is the function field of $X'$,  is hardly ever a neighbourhood of $s$. Nevertheless, we can choose a finite extension $k'/k$ such that $s(\Gal_{k'}) \subset \Gal_{K'}$ which means that a scalar extension to $k'$ (from the field of constants) of $X'$ is a neighbourhood of $s' = s|_{\Gal_{k'}}$ as a section of the birational extension associated to $X \times_k k'$. Now
\[
X_{s'} = X_s \times_k k'
\]
and the map of rational points factors as 
\[
X_s(k') \to X'(k') \to X(k'),
\]
with $a \not=b$ in the image. Let $a'$ and $b'$ be the intermediate images in $X'(k')$. Then we have 
\[
s_{a'} =s_{b'}
\]
as sections of $\pi_1(X'/k')$, so that 
\[
f'(a') = f'(b')
\]
by the injectivity of the map from rational points to sections for abelian varieties over a $p$-adic local field. This contradicts our choice of $f$, namely $f(a) \not= f(b)$, and the proof is complete.
\end{pro}

\medskip

We are now ready to prove the main theorem.

\begin{thm} \label{thm:main}
Let $X/k$ be a geometrically irreducible, normal, proper variety over a finite extension $k/\bQ_p$. Let $K$ be the function field of $X$. Then every  section of the natural projection  
\[
\Gal_K \to \Gal_{k}
\]
has image in the decomposition subgroup $D_{\bar v} \subset \Gal_K$ for a unique $k$-valuation
 $\bar{v}$ of $\ov{K}$ with residue field of $v = \bar{v}|_K$ equal to  $k$. 
 
 In particular, conjugacy classes of sections of $\Gal_K \to \Gal_{k}$ come in disjoint non-empty packets associated to each $k$-valuation $v$ of $K$ with residue field $k$.
\end{thm}
\begin{pro}
We first show that for every $k$-valuation $v$ the set of birational Galois sections with image in $D_v$ is nonempty.
Put differently, we need to show that the natural projection $D_v \to \Gal_k$ splits.  The corresponding extension
\[
1 \to I_v \to D_v \to \Gal_k \to 1
\]
with the inertia group $I_v$ of $v$ splits by general valuation theory. 

\smallskip

We now show that every birational Galois section $s$ belongs to a packet of a $k$-valuation. By Chow's Lemma and Hironaka's resolution of singularities, we may assume that $X/k$ is smooth and projective. For every normal, birational $X' \to X$ the section $s$ gives rise to a tower of branched neighbourhoods that are again linked by a natural map 
\[
X'_s \to X_s.
\]
The map that assigns to a $k$-valuation of $\ov{K}$ its center on a birational model for a finite branched cover of $X$ leads to a bijection
\[
\left\{\bar{v} \ ; \ k\text{-valuation on $\ov{K}$ with residue field algebraic over $k$} \right\} \xrightarrow{\sim} \varprojlim_{X' \to X} X'_s(\bar{k}) 
\]
since $X'_s(\bar{k}) = X_s \times_k \bar{k} (\bar{k})$ and $X'_s \times_k \bar{k}$ is the normalization of $X'_s$ in $\ov{K}$. Indeed, the assumption on the residue field of the valuation implies that the center is a closed point on every birational model. 

Now, a valuation $\bar{v}$ has $s(\Gal_k) \subset D_{\bar{v}}$ if and only if the action of $\Gal_k$ via $s$ on the set of all valuations fixes $\bar{v}$. But the set of such fixed points is precisely
\[
\varprojlim_{X' \to X} X'_s(k) 
\]
which is a set of cardinality one by Proposition \ref{prop:card1}. This proves the theorem.
\end{pro}

\begin{rmk}
I am grateful to Jochen Koenigsmann for allowing me to include his observation that Theorem \ref{thm:main} also admits the following  valuation theoretic proof based on \cite{koenigsmann:p-rigidelements}.

\smallskip

\begin{pro}[Alternate proof] Let $s : \Gal_k \to \Gal_K$ be a section of $\res_{K/k} : \Gal_K \to \Gal_k$, and here $K/k$ can be any extension of fields where $k/\bQ_p$ is finite. Denote again by $L \subseteq \ov{K}$ the fixed field under $s(\Gal_k)$.  Then $\Gal_L \cong \Gal_k$ and, by \cite{koenigsmann:p-rigidelements} Theorem 4.1,
the field $L$ is $p$-adically closed with respect to a unique $p$-adic valuation $w_L$ that extends the unique $p$-adic valuation on $k$. 

Now let $v_L$ be the \textit{relative rank-1} coarsening of $w_L$, i.e., the $k$-valuation with valuation ring $O_{w_L}[1/p]$. 
The valuation $v_L$ is nontrivial, since otherwise $L/k$ would be an  immediate extension for the $p$-adic valuations, and such do not exist nontrivially as $k$ is complete.  As a nontrivial coarsening of a henselian valuation, $v_L$ is henselian as well, say with unique prolongation $\bar v$ to $\ov{K}$ that is a coarsening of the unique prolongation $\bar w$ of $w_F$. We conclude that
\[
s(\Gal_k) = \Gal(\ov{K}/L) \subseteq D_{\bar w} = D_{\bar{v}} 
\]
and it remains to show that $v = \bar{v}|_K$ has residue field $k$. 

The residue field of $L$ with respect to $v_L$ is an algebraic  $p$-adically closed
extension of $k$ with absolute Galois group isomorphic to $\Gal_k$
(note that the value group of $v_L$ is divisible),
hence it is equal to $k$ and this is inherited by the residue field of $v$.
\end{pro}

\medskip

The generalization to arbitrary field extensions could be obtained by a limit argument from Theorem \ref{thm:main} as well.
\end{rmk}

%%%%%%%%%%%%%%%%%%%%%%%%%%%%%%%%%%%%%%%%%%
%%%%%% Bibliography %%%%%%%%%%%%%%%%%%%%%%%%%%%%%
%%%%%%%%%%%%%%%%%%%%%%%%%%%%%%%%%%%%%%%%%%

%%%%%%%%%%%%%%%%%%%%%%%%%%%%%%%%%%%%%%%%%%

\end{document}